\documentclass{article}

\usepackage{arxiv}

\usepackage[utf8]{inputenc} 
\usepackage[T1]{fontenc}    
\usepackage{hyperref}       
\usepackage{url}            
\usepackage{booktabs}       
\usepackage{amsfonts}       
\usepackage{nicefrac}       
\usepackage{microtype}      
\usepackage{lipsum}         
\usepackage{graphicx}
\usepackage{doi}
\usepackage{graphicx}
\usepackage{amssymb}  
\usepackage{amsmath}
\usepackage{pgfplots}
\usepackage{cleveref}       

\title{Optimal Control of Longitudinal Motions for an Elastic Rod with Distributed Forces \thanks{The study has been done under financial support of
	the Russian Science Foundation (grant 21-11-00151).}
}

\author{ \href{https://orcid.org/ 0000-0001-6526-6246}{\includegraphics[scale=0.06]{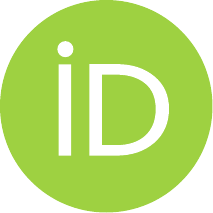}\hspace{1mm}Georgy Kostin}
 \\	Ishlinsky Institute for Problems in Mechanics RAS\\ 
	Moscow, Russia \\	
		\texttt{kostin@ipmnet.ru} \\
	\And
	\href{https://orcid.org/0000-0001-7210-5395}{\includegraphics[scale=0.06]{orcid.pdf}\hspace{1mm}Alexander Gavrikov} \\
	Department of Mathematics\\
	Penn State University\\
	State College, USA \\
	\texttt{avg6113@psu.edu} \\
	Ishlinsky Institute for Problems in Mechanics RAS\\
	Moscow, Russia \\	
}

\renewcommand{\shorttitle}{Optimal Control of Longitudinal Motions}

\begin{document}
\def\ds{\displaystyle}
\maketitle

\begin{abstract}
The study is devoted to mathematical modeling and optimal control design of longitudinal motions of a rectilinear elastic rod. The control inputs are a force, which is normal to the cross section and distributed piecewise constantly along the rod’s axis, as well as two external lumped loads at the ends. It is assumed that the intervals of constancy in the normal force have equal length. Given initial and terminal states with a fixed time horizon, the optimal control problem is to minimize the mean mechanical energy stored in the rod. 
To solve the problem, two unknown functions are introduced: the dynamical potential and the longitudinal displacements. As a result, the initial-boundary value problem is reformulated in a weak form, in which constitutive relations are given as an integral quadratic equation. The unknown functions are both continuous in the new statement. For the uniform rod, they are found as linear combinations of traveling waves. In this case, all conditions on continuity as well as boundary, initial, and terminal constraints form a linear algebraic system with respect to the traveling waves and control functions. The minimal controllability time is found from the solvability condition for this algebraic system.
After resolving the system, remaining free variables are used to optimize the cost functional. Thus, the original control problem is reduced to a one-dimensional variational problem. The Euler–Lagrange necessary condition yields a linear system of ordinary differential equations with constant coefficients supplemented by essential and natural boundary conditions. Therefore, the exact optimal control law and the corresponding dynamic and kinematic fields are found explicitly. Finally, the energy properties of the optimal solution are analyzed.
\end{abstract}


\section{Introduction}

For many years, such classical mechanical systems as thin elastic rods or strings and related control problems have been attracting attention of scientists. 
This is not surprising given how many physical processes are modeled by the wave equation. A possible solution to a control problem for this model involves boundary and distributed control inputs \cite{Lions:1971}.
The boundary control is often more feasible to implement in mechanical systems since its practical realization involves drivers widely used in engineering. However, it has intrinsic limitations since  a finite number of inputs is used to control a continuum system. For a rod, one of the limitations is a minimal control time 
\cite{Butkovsky:1969}, \cite{Kostin:2021}. 


A practical implementation of a distributed control law may suppose preliminary spatial discretization and only then be applied to a mechanical system.
In our study, we assume from the beginning that the control inputs are finite-dimensional: boundary control forces are applied at the ends of the rod and a piecewise constant normal force in the cross section acts along its central line. 
Thus, our distributed control is discrete in space, although each of the inputs can take continuous values in time. Such control structure can be implemented, for example, with a set of piezoelectric actuators placed symmetrically along the entire length of the rod. 
 In this work, we do not consider a detailed model of the actuators and understand forces as control inputs for simplicity. For practical implementation our control strategy may be complemented by a model of piezoelements   as in \cite{Ramesh_Kumar:2008},~\cite{Li:2017}. 

The assumption that the input is piecewise constant allows us to split the system into several interconnected subsystems so that each subsystem may be described by means of auxiliary traveling waves. The continuity and boundary conditions as well as control inputs interweave these functions. We use a special mesh in the time-space domain to express all the conditions in terms of a linear system, which solvability guarantees the controllability of the system. Our approach allows us to derive the exact analytical expressions for optimal control inputs, which is not typical  for distributed parameter systems since usually only a numerical (approximate) solution is possible.

\section{Statement of the Control Problem}\label{sec:2}

\subsection{Rod Kinematics and Control Forces}\label{sub:201}

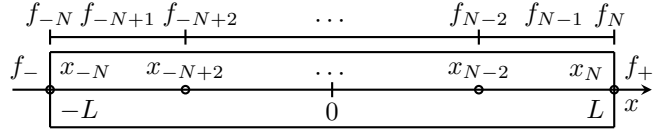
\begin{figure}[t]
	\begin{center}
		\begin{tikzpicture}[thick,scale=1.0, every node/.style={transform shape}]
			\draw [->,>=stealth] (0,0) -- (8.5,0);
			\node[below right] at (8,0){$x$}; 
			\draw (0.5,0.5) -- (8,0.5);
			\draw (0.5,-0.5) -- (8,-0.5);
			\draw (0.5,-0.5) -- (0.5,0.5);
			\draw (8,-0.5) -- (8,0.5);
			\node[below right] at (0.5,0){$-L$}; 
			\node[below left] at (8,0){$L$}; 
			\node[above right] at (0.5,0){$x_{-N}$};
			\node[above left] at (8,0){$x_{N}$}; 
			\node[above] at (2.3,0){$x_{-N+2}$}; 
			\node[above] at (6.2,0){$x_{N-2}$}; 
			\draw (4.25,-0.1) -- (4.25,0.1);
			\draw (0.5,0) circle (.05 cm);
			\draw (8,0) circle (.05 cm);
			\draw (2.3,0) circle (.05 cm);
			\draw (6.2,0) circle (.05 cm);
			\node[below] at (4.25,-0.05){$0$}; 
			\node[above] at (4.25,0){$\cdots$}; 
			\draw (0.5,0.7) -- (8,0.7);
			\draw (0.5,0.6) -- (0.5,0.8);
			\draw (8,0.6) -- (8,0.8);
			\draw (2.3,0.6) -- (2.3,0.8);
			\draw (6.2,0.6) -- (6.2,0.8);
			\node[above] at (0.5,0.7){$f_{-N}$}; 
			\node[above] at (1.4,0.7){$f_{-N+1}$}; 
			\node[above] at (2.5,0.7){$f_{-N+2}$}; 
			\node[above] at ((4.25,0.7){$\cdots$};
			\node[above] at (6.2,0.7){$f_{N-2}$}; 
			\node[above] at (7.2,0.7){$f_{N-1}$}; 
			\node[above] at (7.95,0.7){$f_{N}$}; 
			\node[above right] at (8,0){$f_{+}$}; 
			\node[above left] at (0.5,0){$f_{-}$}; 
		\end{tikzpicture}\\
	\end{center}  
	\caption{Scheme of a rod with control elements.}
	\label{fig:01}
\end{figure}

Let us consider longitudinal motions
of a uniform rectilinear elastic rod of length $2L$ (see Fig.~\ref{fig:01}). The $x$-axis is directed along the central line with the origin at the middle of the rod.
The absolute displacement of a point with the coordinate $x\in I_L$, $I_L= (-L,L)$ at the time instant $t\in I_T$, $I_T=(0,T)$ is given by a mapping $v:\Omega\rightarrow\mathbb{R}$, where $\Omega=I_T\times I_L$ is the time-space domain.     

The rod is controlled by two external normal forces $f_{\pm}:I_T\rightarrow\mathbb{R}$ applied at the ends with $x=\pm L$.
Additionally, the force $f:\Omega\rightarrow\mathbb{R}$ normal to the cross section tenses or compresses the elastic rod along the $x$-axis. We assume that  the function $f(t,x)$ is piecewise constant in space. 
Such a load can be generated, for example, by a set of piezoelectrical actuators (control elements) attached  to the rod's lateral surface \cite{Ramesh_Kumar:2008}, \cite{Li:2017}. 
In our setting, these $N$ elements have equal length and are inseparably located along the central line.

The control elements  are naturally related to $N$ space intervals $I^{x}_{k}$ defined along the rod by  points  $x_{n}$ given by
\begin{equation} \label{eq:203:rod_segments}
	\begin{array}{c}
		x \in I^{x}_{k} := (x_{k-1},x_{k+1}),\quad
		k \in J_s;
	\\ 
		x_{n} =n\lambda/2,\quad
		n \in J_x,\quad 
		\lambda = 2L/N,\quad 
		x_{\pm N} = \pm L,
	\end{array}		
\end{equation} 
where $\lambda$ is the length of each segment.
The  sets of indices 
\begin{equation} \label{eq:203:space_indexing_sets}
	\begin{array}{c}
		J_{s} = \{ 1-N, 3-N,\ldots,N-1 \},
	\quad 
		J_{x} = \{ -N, 2-N,\ldots,N \}
	\end{array}		
\end{equation} 
are introduced in~\eqref{eq:203:rod_segments} to specify respectively the space intervals $I^{x}_{k}$ and the interface points $x_{n}$. 

The control function $f$ over each segment $x \in I^{x}_{k}$ is expressed in the form $f(t,x) = f_{k}(t)$ with $k\in J_s$, $t\in  I_T$.
The renamed external loads $f_{-N-1}:=f_{-}$ and $f_{N+1}:=f_{+}$ complete a set of control functions $f_{k}$ with the indices  $k\in J_c= J_{s}\cup\{-N-1, N+1\}$. 
Here the supplemented index set is related to the set of control inputs (see Fig.~\ref{fig:01}).

\subsection{Initial-Boundary Value Problem}\label{sub:204}

According to the method of integro-differential relations (MIDR), \cite{Kostin:2018}, we define the linear momentum density $p:\Omega\rightarrow\mathbb{R}$ and the total internal normal forces in the cross section $s:\Omega\rightarrow\mathbb{R}$ as unknown variables together with the displacements $v$.
The governing PDEs 
can be split into two parts. 
The first part is a PDE, which links the momentum and the force 
\begin{equation} \label{eq:204:Newton_law}
	p_{t}(t,x) = s_{x}(t,x),\quad
	(t,x) \in \Omega
\end{equation}
according to Newton's second law. 
The subscripts $t$ and $x$ denote the partial derivatives in time and space. 
The second part includes two constitutive relations
\begin{equation} \label{eq:204:constitutive_law}
	\begin{array}{l}
	g(t,x) = 0, \quad
	h(t,x) = 0, \quad 
	(t,x)\in \Omega,
	\\
	g := \rho v_{t} - p, \quad
	h := \kappa v_{x} - s + f,
\end{array}
\end{equation}
between the momentum $p$ and the velocities $v_t$ as well as
between the forces $s$ and the longitudinal strains $v_x$.
This model contains  two constants, namely, 
the linear density $\rho>0$ and the tension stiffness $\kappa>0$ of the rod.

The initial conditions are imposed
on both the displacements and the momentum density by
\begin{equation} \label{eq:204:initial_distributions}
	v(0,x) = v_{0}(x),\quad
	p(0,x) = p_{0}(x), \quad
	x \in  I_L,
\end{equation}
where $v_0(x)$ and $p_0(x)$ are given functions.   
The boundary conditions are defined as follows
\begin{equation} \label{eq:204:boundary_distributions}
	s(t, -L) = f_{-N-1}(t),\;\;
	s(t,  L) = f_{N+1}(t), \;\;
	t \in  I_T.
\end{equation}
The continuity conditions at the inner interface points
\begin{equation} \label{eq:204:interelement_distributions}
	\begin{array}{l}
		[v(t,x_{n})] = 0, \quad
		[s(t,x_{n})] = 0, 
	\quad 
		n\in J_{x} \setminus \{ -N, N \}, \quad
		t \in  I_T,
	\end{array}
\end{equation}
for displacements and forces must be also taken into account. 
Here the square brackets denote the jump value of the variable $v$ or $s$. 
For shortness, we define the relations~\eqref{eq:204:initial_distributions}--\eqref{eq:204:interelement_distributions} as interface conditions.

The jumps of derivatives $[v_x(t,x_{n})]$ and the boundary values of these derivatives $v_x(t,x_{\pm N})$ depend exclusively on the differences of adjacent control functions
\begin{equation} \label{eq:204:jump_force}
	f_{n}:=f_{n+1}-f_{n-1},\quad
	n \in J_x.
\end{equation}
The number of these function is one less than the number of the original inputs $f_{k}$ with $k \in J_c$. 
This means that the dynamic behavior of the system does not depend on the sum of the control signals $f_{k}(t)$. 
Its value affects only the intensity of residual stresses in the rod material. 
To reduce residual stresses, let us set this sum to zero as follows
\begin{equation} \label{eq:204:zero_force_sum}
	\sum\nolimits_{k \in J_c} f_{k}(t) = 0, \quad
	t\in I_T.
\end{equation}
We define two vector spaces with elements $\boldsymbol{f}_{c}: I_T \rightarrow \mathbb{R}^{N+2}$ (boundary and distributed control inputs) and $\boldsymbol{f}:  I_T \rightarrow \mathbb{R}^{N+1}$ (control jumps) described according to
\begin{equation} \label{eq:204:force_vectors}
	\begin{array}{c}
		\boldsymbol{f}_c := (f_{-N-1}, f_{-N+1},\ldots, f_{N-1}, f_{N+1}) , 
		\\ 
		\boldsymbol{f} := (f_{-N}, f_{-N+2},\ldots, f_{N-2}, f_{N}) .
	\end{array}
\end{equation}
Given the vector $\boldsymbol{f}$, we can uniquely resolve the linear system of 
 equations~\eqref{eq:204:jump_force}, \eqref{eq:204:zero_force_sum} w.r.t. the entries of $\boldsymbol{f}_{c}$. 

\subsection{Dynamic Potential in an Equivalent IBVP Statement}\label{sub:205}

Next, we  define an auxiliary function or a   dynamic `potential'  $r:\Omega \rightarrow \mathbb{R}$ such that
\begin{equation} \label{eq:205:potential_relations}
	p(t,x) = r_x(t,x), \quad
	s(t,x) = r_t(t,x), \quad 
	(t,x)\in \Omega.
\end{equation}
This representation of the momentum density $p$ and the normal forces $s$ automatically satisfies~\eqref{eq:204:Newton_law}. 
By excluding $p$ and $s$  from consideration, the state of the system is determined by two variables, kinematic $v$ and dynamic $r$.

To complete properly the PDE system~\eqref{eq:204:constitutive_law}, we have to rewrite the interface conditions~\eqref{eq:204:initial_distributions}--\eqref{eq:204:interelement_distributions}	 in terms of the new variables $v$ and $r$. 
As only the partial derivatives $r_t$ and $r_x$ appear in the governing equations~\eqref{eq:204:constitutive_law}--\eqref{eq:204:interelement_distributions},
the variable $r(t, x)$ can be found up to an arbitrary constant $c_{0,k}$ on each subdomain $\Omega_{k} =I_T\times I^{x}_{k} \subset \Omega$ with $k \in J_s$.

By imposing continuity conditions $[r(0,x_{n})] = 0$ with $n\in J_{x} \setminus \{ -N, N \}$ on the initial values of dynamic variable $r$, the number of constants can be reduced to one ($c_0$)  and the initial conditions~\eqref{eq:204:initial_distributions} are rewritten as follows
\begin{equation} \label{eq:205:initial_conditions}
		v(0,x) = v_{0}(x),\quad
		r(0,x) = r_{0}(x)
		=\int_{-L}^{x}p_0(\xi)\,\mathrm{d}\xi + c_{0}.
\end{equation}
After integrating~\eqref{eq:204:boundary_distributions} and~\eqref{eq:204:interelement_distributions} in view of~\eqref{eq:205:potential_relations}, the new boundary constraints for $t \in I_T$ have the form 
\begin{equation} \label{eq:205:boundary_conditions}
	r(t,\pm L) = r_{0}(\pm L) + u_{\pm N \pm 1}(t).
\end{equation}	
The new interelement conditions are
\begin{equation} \label{eq:205:interelement_conditions}
	[v(t,x_{n})] = 0, \quad
	[r(t,x_{n})] = 0, \quad
	n \in J_{x} \setminus \{ -N, N \}.
\end{equation}

The pair $u_{\pm N\pm 1}$ in~\eqref{eq:205:boundary_conditions} is control integrals given by 
\begin{equation} \label{eq:206:control_integrals}
	u_{k}(t)=\int_{0}^{t} f_{k}(\tau)\,\mathrm{d}\tau, \quad
	k \in J_c,\quad
	t \in I_T.
\end{equation}
Similarly to~\eqref{eq:204:jump_force}, 
 we can also define the jumps of control integrals $u_{n}(t)=u_{n+1}(t)-u_{n-1}(t)$ with $n \in J_x$.
According to~\eqref{eq:204:force_vectors}, \eqref{eq:206:control_integrals}, two vector-valued functions $\boldsymbol{u}_{c}:I_T\rightarrow \mathbb{R}^{N+2}$ and $\boldsymbol{u}:I_T\rightarrow \mathbb{R}^{N+1}$ are introduced.

\subsection{Variational Formulation and Optimal Control Problem}\label{sub:207}

Following the MIDR, an integral relation replaces constitutive PDEs to derive a variational statement in dynamics.
The variational reformulation of the dynamic problem~\eqref{eq:204:constitutive_law},  \eqref{eq:205:initial_conditions}--\eqref{eq:205:interelement_conditions} can be described as follows.
Given the values of control function $\boldsymbol{u}(t)$ for $t\in I_T$, find such functions $v^{*}(t,x,\boldsymbol{u})$ and $r^{*}(t,x,\boldsymbol{u})$ in the Hilbert space $H^{1}(\Omega)$
that minimize the value of the constitutive functional
\begin{equation} \label{eq:207:constitutive_minimization}
	\begin{array}{c}
\displaystyle		Q[v^{*},r^{*},\boldsymbol{u}] 
		= \min\limits_{v,r\in H^1(\Omega)}Q[v,r,\boldsymbol{u}]=0,
	\quad
		Q = \int\nolimits_{\Omega}q\,\mathrm{d}\Omega,\quad
		q=\frac{ g^{2}}{4\rho}+\frac{h^{2}}{4\kappa}, 
	\end{array}
\end{equation}
subject to the interface constraints \eqref{eq:205:initial_conditions}--\eqref{eq:205:interelement_conditions}. 
Here the constitutive residual function $q$ follows from~\eqref{eq:204:constitutive_law} and~\eqref{eq:205:potential_relations}. 
We consider  $\boldsymbol{u}$ as a frozen parameter  in~\eqref{eq:207:constitutive_minimization}.   

The following optimal control problem is considered.
Find control vector-valued function $\boldsymbol{u}^{*}(t)$ such that the mean mechanical energy $E$  stored in the rod over the fixed time horizon $ I_T$ reaches its minimum: 
\begin{equation} \label{eq:209:energy_minimization}
	\begin{array}{c}
\displaystyle		E[v,r,\boldsymbol{u}]\to \min\limits_{\boldsymbol{u},c_{1}},\quad
		E = \frac{1}{T}\int\nolimits_{\Omega}e \,\mathrm{d}\Omega, 
		\quad 
		e=\frac{\rho v_{t}^{2}}{4}+\frac{\kappa v_{x}^{2}}{4}
		+\frac{(r_{t}-f)^{2}}{4\kappa}+\frac{r_{x}^{2}}{4\rho}, 
	\end{array}
\end{equation}
subject to the integral constraint~\eqref{eq:207:constitutive_minimization},
the interface conditions~\eqref{eq:205:initial_conditions}--\eqref{eq:205:interelement_conditions}, 
and the continuous terminal conditions
\begin{equation} \label{eq:209:terminal_conditions}
		v(T,x)=v_{1}(x),\quad
		r(T,x)=r_{1}(x)
		=\int_{-L}^{x}p_{1}(\xi)\mathrm{d}\xi+c_{1},\quad x\in[-L,L].
\end{equation}
The mapping $e:\Omega \rightarrow \mathbb{R}$ denotes the linear energy density.
The desired functions of displacements $v_{1}:[-L,L]\rightarrow\Omega$ and momentum density $p_{1}:[-L,L]\rightarrow\Omega$ completely define the terminal state of the  rod, whereas the parameter $c_{1}$ does not influence on this state.
The constant $c_{0}$  in~\eqref{eq:205:initial_conditions} is another free parameter. 
A linear combination of $c_{0}$ and $c_{1}$ must be a variable parameter of optimization. 
Without lost of generality, we can  set $c_{0} = 0$.

\section{Solution of the Direct Problem}\label{sec:3}
\subsection{Representation of the Solution in d'Alembert's Form}\label{sub:301}

To simplify the description of the system, we use dimensionless variables according to	
$v= v^{*}$, $r=\kappa \tau_{*}r^{*}$, $x=L x^{*}$, $t=\tau_{*} t^{*}$, $\tau^2_{*}=L^2\rho/\kappa$. 
The star superscript is further omitted. After this transformation, the length of the rod is equal to 2, whereas the length of each segment is $\lambda = 2/N$.

To analyze the optimal dynamics of a rod described in the previous section, a traveling wave representation of the problem variables $(v, r)$ in d'Alembert's form is utilized. 
On each subdomain $\Omega_{k}= I_T\times I^{x}_{k} \subset \Omega$ with $k \in J_{s}$, the kinematic and dynamic variables are represented as
\begin{equation} \label{eq:302:dAlembert_solutions}
	\left\{ 
	\begin{array}{l}
		v(t,x) = w^{+}_{k}(t+x) + w^{-}_{k}(t-x)
		\\ 
		r(t,x) = w^{+}_{k}(t+x) - w^{-}_{k}(t-x)  - u_{k}(t)
	\end{array}
	\right.,
\end{equation}  
where left and right traveling wave functions
$w^{\pm}_{k}: I^{\pm}_{k} \rightarrow \mathbb{R}$, $w^{\pm}_{k}\in C(I^{\pm}_{k})$, 
are introduced with the domains
\begin{equation} \label{eq:302:dAlembert_domains} 
	\begin{array}{l}\textstyle
		I^{\pm}_{k} = \left( z^{\pm}_{k}, T-z^{\mp}_{k} \right),\quad
		z^{+}_{k} = \frac{k-1}{2}\lambda,\quad
		z^{-}_{k} = -\frac{k+1}{2}\lambda.
	\end{array}
\end{equation}

\begin{figure}[t]
	\begin{center}
		\includegraphics[width=0.6\linewidth]{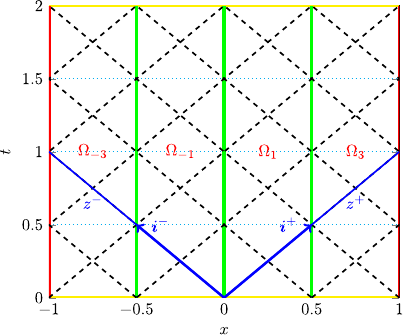}\\ 
		\caption{A mesh on the dimensionless time-space domain $\Omega = (0, 2)\times(-1,1)$ with $M=4$, $N=4$.}
		\label{fig:02}
	\end{center}
\end{figure}

A rather conventional coordinate representation of the traveling waves $w^{\pm}_{k}$ for $k\in J_s$  is given in the coordinate frame $(z^{+}, z^{-})$ rotated counter-clockwise on the angle $\frac{\pi}{4}$ w.r.t. the frame $(t, x)$.
The transformation from the old coordinates to new ones has the form
		$z^{+} = t + x$, $z^{-} = t - x$.
In Fig.~\ref{fig:02}, the $z^{\pm}$-axes are presented by blue lines, the orts of the new frame
$\boldsymbol{i}^{\pm}$ with the coordinates $ <\frac{1}{2},\pm\frac{1}{2}>$ in the old system $(t,x)$ are depicted by two blue arrows. 
In the new coordinates, the traveling functions $w^{\pm}_{k}$ with $k \in J_s$ depend respectively on the only argument $z^{\pm}$.

Although the variables $(v,r)$ in the form~\eqref{eq:302:dAlembert_solutions} satisfy the integral constraint $Q=0$, the solution has been given so far only for the disjoint open subdomains $\bigcup_{k \in J_{s}} \Omega_{k}$. 
By taking into account the initial and terminal constraints ~\eqref{eq:205:initial_conditions}, \eqref{eq:209:terminal_conditions} as well as the boundary and interelement constraints~\eqref{eq:205:boundary_conditions}, \eqref{eq:205:interelement_conditions}, the solution $(v,r)$ needs to be extended to the set of interface edges $\overline{\Omega} \setminus \bigcup_{k \in J_s} \Omega_{k}$ presented in Fig.~\ref{fig:02} for $N=4$ with yellow, red, and green lines.

To resolve all the above-mentioned continuity conditions, we consider only the control time horizon $T$ which is a multiple of the segment length $\lambda$, that is $T=M\lambda$ with $M\in \mathbb{N}$. 
The solution for the other values of $T$ is found similarly but requires more complicated algorithm and is out the scope of the paper.   
The mesh on the time-space rectangle $\overline{\Omega}$ is generated by the interface line segments $t=0$, $t=T$, and $x=x_{n}$ for $n\in J_x$ together with the characteristic segments $z^{\pm}=\frac{i}{2}\lambda$ for $i = N+2j$ with $j \in \mathbb{Z}$.
In Fig.~\ref{fig:02}, the characteristic segments of the mesh with $M=N=4$ are depicted by dashed lines. 
At that, $\lambda=\frac12$ and $T=4\lambda=2$ for this domain.

Similarly to the space intervals $I^{x}_{k}$ and the points $x_{n}$ introduced in~\eqref{eq:203:rod_segments}, the duration intervals $I^{t}_{l}$ and the time instants $t_{m}$ generated by the mesh are given by 
\begin{equation} \label{eq:305:time_intervals}
	\begin{array}{c}
		t \in I^{t}_{l} := (t_{l-1},t_{l+1}),\quad 
		l \in J_{d};
		\\ \textstyle 
		t_{m} = \frac{1}{2}m\lambda,\quad 
		m \in J_{t},  \quad 
		\lambda = 2/N, \quad
		t_{2N} = T.
	\end{array}		
\end{equation} 
The two new sets of indices  are introduced in~\eqref{eq:305:time_intervals}
\begin{equation} \label{eq:305:t_indexing_sets}
	\begin{array}{c}
		J_{d} := \{ 1, 3,\ldots, 2M-3, 2M-1 \},
		\quad 
		J_{t} := \left\{ 0, 2,\ldots, 2M-2, 2M \right\}.
	\end{array}		
\end{equation} 

\subsection{Double Indexing of D'Alembert's and Control Functions}\label{sub:306}

To extend the solution $(v,r)$ over each closed set $\overline\Omega_{k}$, where $\Omega_{k} =  I_T\times I^{x}_{k}$ with $k\in J_s$, we need to satisfy the continuity conditions over the initial and terminal edges 
\begin{equation} \label{eq:306:time_edges}
	\left\{
	\begin{array}{l}
		S^{0}_{k} = \left\{(t,x)\in\overline\Omega_k: t=0 \right\},
		\\ 
		S^{T}_{k} = \left\{(t,x)\in\overline\Omega_k: t=T \right\},
	\end{array}		
	\right.\quad
	k \in J_s,
\end{equation} 
and the edges parallel to the $t$-axis 	
\begin{equation} \label{eq:306:space_edges}
	\begin{array}{c}
		S^{x}_{n,l} = \left\{ (t,x)\in\overline\Omega: t\in I^{t}_{l}, x = x_{n} \right\},\quad
		n\in J_{x},\quad
		l\in J_{d}.
	\end{array}		
\end{equation} 		

To operate with the values of $w^{\pm}_{k}$ on these edges in $\overline{\Omega}_{k}$, the open intervals $I^{\pm}_{k,m}$ of the length $\lambda$ are selected on each domain $I^{\pm}_{k}$ from~\eqref{eq:302:dAlembert_domains} according to
\begin{equation} \label{eq:306:wave_function_intervals}
	\begin{array}{c}
		I^{\pm}_{k,m} = ( z^{\pm}_{k,m},  z^{\pm}_{k,m+2} ),\quad
		z^{\pm}_{k,m} = z^{\pm}_{k} + \frac{m}{2}\lambda,\;
		m\in J_{t}.
	\end{array}		
\end{equation} 		
Here, $J_t$ is the set of indices in~\eqref{eq:305:t_indexing_sets}, and the characteristic coordinate $z^{\pm}_{k}$ is given in~\eqref{eq:302:dAlembert_domains}.	

The new functions $w^{\pm}_{k,m}:(0,\lambda)\rightarrow \mathbb{R}$ are defined so that
\begin{equation} \label{eq:306:edge_wave_functions}
	\begin{array}{c}
		w^{\pm}_{k,m}(z) = w^{\pm}_{k}( z+z^{\pm}_{k,m} ),\quad
		k \in J_{s},\quad
		m \in J_{t}, 
	\end{array}		
\end{equation} 		
where the coordinate shifts $z^{\pm}_{k,m}$ are represented in~\eqref{eq:306:wave_function_intervals}, the index sets $J_{s}$, $J_{t}$ are introduced in~\eqref{eq:203:space_indexing_sets}, \eqref{eq:305:t_indexing_sets}, respectively. 
The same procedure is applied to the control functions $u_{j}(t)$ given on $t\in[0,T]$. 
Each of them is split into $m$ edge maps $u_{j,m}:(0,\lambda)\rightarrow \mathbb{R}$, which values are equal to
\begin{equation} \label{eq:306:edge_control_functions}
	\begin{array}{c}
		u_{j,m}(z) = u_{j}\left( z + t_{m} \right),\quad
		j\in J_{x}\cup J_{c},\quad
		m \in J_{t}  \setminus \{ 2M \},
	\end{array}		
\end{equation} 	
where the time instants $t_{m}$ are given in~\eqref{eq:305:time_intervals}.	

After accounting for the zero initial values of $u_{k}$, the initial conditions from~\eqref{eq:205:initial_conditions} in d'Alembert's form~\eqref{eq:302:dAlembert_solutions} are given on the mesh edges $S^{0}_{k}$ introduced in~\eqref{eq:306:time_edges} by   
\begin{equation} \label{eq:307:edge_initial_conditions} 
	\left\{
	\begin{array}{l}
		w^{+}_{k,0}(z) + w^{-}_{k,0}(\lambda-z) = v_{0}(z^{+}_{k}+z) 
		\\ 
		w^{+}_{k,0}(z) - w^{-}_{k,0}(\lambda-z) = r_{0}(z^{+}_{k}+z)
	\end{array}\right.,\quad
	k\in J_s.  
\end{equation}  
For the terminal conditions~\eqref{eq:209:terminal_conditions} on the edges $S^{T}_{k}$, we similarly arrive at	
\begin{equation} \label{eq:307:edge_terminal_conditions} 
	\left\{
	\begin{array}{l}
		w^{+}_{k,2M}(z) + w^{-}_{k,2M}(\lambda-z) = v_{1}(z^{+}_{k}+z) 
		\\ 
		w^{+}_{k,2M}(z) - w^{-}_{k,2M}(\lambda-z) = r_{1}(z^{+}_{k}+z)
	\end{array}\right.,\quad
	k\in J_s.  
\end{equation}  
Note that \eqref{eq:307:edge_initial_conditions} and~\eqref{eq:307:edge_terminal_conditions} contains $4N$ relations. 

The boundary conditions on $S^{x}_{\pm N, m}$ from~\eqref{eq:306:space_edges} bind the traveling waves and the boundary control jumps as follows
\begin{equation} \label{eq:307:edge_boundary_conditions} 
	\begin{array}{c}
		w^{-}_{1-N,m+2}(z)-w^{+}_{1-N,m}(z) = u_{-N,m}(z)+r_{0}(-1), 
	\\ 
		w^{+}_{N-1,m+2}(z) - w^{-}_{N-1,m}(z) = u_{N,m}(z)+r_{0}(1)
	\end{array}
\end{equation}  
with $m\in J_t \setminus \{2M\}$. 
The number of the boundary relations is equal to $2M$.
The continuity conditions on the other segments $S^{x}_{n, m}$ link the wave functions with the control jumps according to 
\begin{equation} \label{eq:307:edge_interelement_conditions} 
	\begin{array}{c}
			w^{+}_{n-1,m+2} + w^{-}_{n-1,m} =
			w^{+}_{n+1,m} + w^{-}_{n+1,m+2}, 
		\\ 
			w^{+}_{n-1,m+2} - w^{-}_{n-1,m} =
			w^{+}_{n+1,m} - w^{-}_{n+1,m+2}  + u_{n,m},
		\\ 
		n\in J_{x} \setminus \{-N,N\},\quad
		m\in J_{t} \setminus \{2M\}.  
	\end{array}
\end{equation}  
There are $2M(N-1)$ equations related to these edges.

Thus, the total number of the edge constraints is equal to $N_{e} = 2MN + 4N$. 
The equations~\eqref{eq:307:edge_initial_conditions}--\eqref{eq:307:edge_interelement_conditions} 
represent a linear algebraic system w.r.t. unknown wave ($w^{\pm}_{i,j}$) and control ($u_{k,l}$) functions. 
The total number of the unknown variables equals to $N_{v} = N_{w} + N_{u} = 3MN + M + 2N$ with the number of the wave functions $N_{w} = 2(M+1)N$ and the number of the control functions  $N_{u} = M(N+1)$. 
As follows from these quantities for the system~\eqref{eq:307:edge_initial_conditions}--\eqref{eq:307:edge_interelement_conditions}, the number of variables exceeds the number of equations for a rather large value of $MN$.  The issues of solvability of this system is discussed in the rest of this section.

\subsection{Solvability of the System of Edge Constraints}\label{sub:308}

The forced motions of an elastic rod with one interval of constancy of the distributed force, that is $N=1$, are equivalent to motions of a rod controlled by external boundary forces $f^{\pm}$ and was described in \cite{Kostin:2021}. 
Here, we suppose  that $N>1$ and $M>1$.
The variable surplus is equal to 
\begin{equation} \label{eq:308:uniform_variable_surplus} 
	N_{s}:=N_{v} - N_{e},\quad 
	N_{s}(M,N)= MN + M - 2N.
\end{equation}  
The function $N_{s}$ is positive and monotonically increasing w.r.t. 
 $M$ if the number of space elements $N$ is fixed. Thus, the  system~\eqref{eq:307:edge_initial_conditions}--\eqref{eq:307:edge_interelement_conditions} is underdetermined.


To prove that the underdetermined inhomogeneous linear system~\eqref{eq:307:edge_initial_conditions}--\eqref{eq:307:edge_interelement_conditions} for $M>1$ and $N>1$ can be resolved, we consider the following scheme. 
Assume that $N=2j+1$ with $j\in \mathbb{N}$, then: 
\begin{enumerate}
\item  the initial and terminal conditions~\eqref{eq:307:edge_initial_conditions}, \eqref{eq:307:edge_terminal_conditions} are resolved  w.r.t. $2N$ functions $w^{\pm}_{k,0}$ and $w^{\pm}_{k,2M}$ for $k\in J_s$;  
\item  the control functions $u_{\pm N,m}$ are found through the boundary conditions~\eqref{eq:307:edge_boundary_conditions} for $m\in J_t \setminus \{2M\}$; \item the interelement conditions~\eqref{eq:307:edge_interelement_conditions} are resolved on 
the edges $S^{x}_{\pm n,m}$ with $m=1,2M-1$ and $n> 0$, namely, the first equation in~\eqref{eq:307:edge_interelement_conditions} is solved w.r.t. $w^{\mp}_{\pm n \pm 1,2}$, $w^{\pm}_{\pm n \pm 1,2M-2}$, while the second equation is solved w.r.t. the control function $u_{\pm n,l}$; 
\item the equations~\eqref{eq:307:edge_interelement_conditions} are resolved on the edges $S^{x}_{\pm n,m}$ with $m\ne 1,2M-1$ and $n> 0$ w.r.t. $w^{\pm}_{\pm n \pm 1,m}$, $w^{\mp}_{\pm n \pm 1,m+2}$; 
\item if $N$ is even, additionally the equations~\eqref{eq:307:edge_interelement_conditions} are resolved  on the central line $S^{x}_{0}$ 
w.r.t. $w^{+}_{-1,0}$, $w^{-}_{1,0}$, $w^{+}_{1,m}$ and $w^{-}_{-1,m}$ with $m\in J_{t} \setminus \{0,2M\}$. 
\end{enumerate}

The proposed scheme for the linear equations~\eqref{eq:307:edge_initial_conditions}--\eqref{eq:307:edge_interelement_conditions} confirms their solvability for the mesh parameters $N>1$ and $M>1$. 
The unresolved variables are combined in a vector $\boldsymbol{y}:(0,\lambda)\rightarrow \mathbb{R}^{N_{s}}$. 

The lack of variables is always present in the case $M=1$ as follows from expression for the variable surplus $N_s$:  $N_{s} = MN + M - 2N$. 
Indeed,  equations \eqref{eq:307:edge_interelement_conditions} will include variables $w^{+}_{n-1,2}$, $w^{-}_{n-1,0}$, $w^{+}_{n+1,0}$, $w^{-}_{n+1,2}$ with $n\in J_{x} \setminus \{-N,N\}$,
which are expressed through the initial and terminal functions $v_{0}$, $r_{0}$, $v_{1}$, $r_{1}$.  
It is always possible to find such a time-boundary distribution of the rod kinematic and dynamic variables $v$ and $r$ that some of these constraints are violated.	
This confirms impossibility to bring the system from an arbitrary initial state to a desired state in the critical time $T=\lambda$.

\subsection{Continuity Conditions at Mesh Vertices}\label{sub:310}

Additionally to the continuity conditions~\eqref{eq:307:edge_initial_conditions}--\eqref{eq:307:edge_interelement_conditions} on  mesh edges, the corresponding linking of governing functions should be done at mesh vertices.
It suffices to continuously conjugate only the free functions with their neighbors. 
If $N$ is odd then the vertex conditions are given by
\begin{equation} \label{eq:310:odd_links}
	\begin{array}{c}
		u_{n,m}(0)=u_{n,m-2}(\lambda),\quad w^{\pm}_{0,m}(0) = w^{\pm}_{0,m-2}(\lambda),\\
		n \in J_{x} \setminus \{ -N, N \}, \quad
		m \in J_{t} \setminus \{0, 2M \}. 
	\end{array}
\end{equation}  
Thus, the total number of equations is equal to $N_{b} = MN+M-N+1$. 
For the even number $n$, the following conditions are added
\begin{equation} \label{eq:310:even_links} 
		u_{0,0}(0)=0, \quad
		u_{0,2M}(0)=u_{0,2M-2}(\lambda); \;  
\end{equation}  
For this case, the number of equations is equal to $N_{b} = MN + M - N$. It should be noted that the equations of the system~\eqref{eq:310:odd_links}, \eqref{eq:310:even_links} are not all necessarily linearly independent 
because some free functions may be coupled with dependent ones, which are, in turn, expressed through the other free variables.

\section{Optimal Control Design}\label{sec:4}

\subsection{Decomposition of the Mean Energy}\label{sub:401}

Noticing that $r_t - f = v_x$ and $r_x = v_t$ on the solution, we obtain that the energy density function $e$ defined in~\eqref{eq:209:energy_minimization} is expressed as $e = \frac12 v_{t}^2 +\frac12 v_{x}^2$ and is given on each  $\Omega_{k}$  by	
\begin{equation}\label{eq:401:energy_density_decomposition}
		e = 		(w^{+\prime}_{k}(z^{+}))^2+((w^{-\prime}_{k}(z^{-}))^2=e^{+}_{k}(z^{+}) +e^{-}_{k}(z^{-}).       
\end{equation}	
It follows that the functional $E$ of the rod mean mechanical energy can be split, in turn, 
into independent parts 	
\begin{equation}\label{eq:401:energy_decomposition}
	\begin{array}{c}
\displaystyle		E =	\sum\limits_{k\in J_s}\big(E^{+}_{k}+ E^{-}_{k}\big),\quad
		E^{\pm}_{k} = \frac{1}{T}\int\nolimits_{\Omega_{k}}  e^{\pm}_{k}(z^{\pm}) \,\mathrm{d}\Omega.
	\end{array}
\end{equation}	
Here, the functional $E^{\pm}_{k}$ is defined if the wave function $w^{\pm}_{k}$ is given over its domain $I^{\pm}_{k}\in \mathbb{R}$ from~\eqref{eq:302:dAlembert_domains}, the set of segment indices $J_s$ is introduced in~\eqref{eq:203:space_indexing_sets}. 	
Substituting the expression for $e^{\pm}_{k}$ from~\eqref {eq:401:energy_density_decomposition} into $ E^{\pm}_{i}$ we arrive to 
\begin{equation}\label{eq:401:term_mean_energy} 
	\begin{array}{l}
\displaystyle		E^{\pm}_{k} = 
		\frac{1}{T}\int\nolimits_{z^{\pm}_{k}}^{T-z^{\mp}_{k}} \left(w^{\pm\prime}_{k}(z^{\pm})\right)^2 \Delta z^{\mp}_{k}(z^{\pm}) \,\mathrm{d}z^{\pm},
	\end{array}
\end{equation}	
where the functions $\Delta z^{\pm}_{i}$, piecewise linear in $z^{\pm}$, arise as a result of primary integration over the coordinate $z^{\mp}$.

Dividing the intervals of integration $I^{\pm}_{k}=[z^{\pm}_{k},T-z^{\mp}_{k}]$ in~\eqref{eq:401:energy_density_decomposition} into $M+1$ equal parts and replacing the function $w^{\pm}_{k}$ for $w^{\pm}_{k,m}$ in accordance with~\eqref{eq:306:wave_function_intervals} on the subintervals $I^{\pm}_{k,k}$ from~\eqref{eq:306:wave_function_intervals}, we arrive to 
\begin{equation}\label{eq:402:energy_decomposition}
	\begin{array}{c}
\displaystyle		E=
		\frac{1}{T}\int\nolimits_{0}^{\lambda} \big(\boldsymbol{G}(z) \boldsymbol{w}^{\prime}(z)\big)\cdot \big(\boldsymbol{G}(z) \boldsymbol{w}^{\prime}(z)\big) dz.
	\end{array}
\end{equation}			
Here, $\boldsymbol{G}:[0,\lambda]\to \mathbb{R}^{N_{w}\times N_{w}}$ is some diagonal weighting matrix-valued function with $N_{w}=2(M+1)N$.
The vector-valued function $\boldsymbol{w}:[0,\lambda]\to \mathbb{R}^{N_{w}}$ is introduced with its elements 
$w_{j}$ given as follows 
\begin{equation}\label{eq:402:traveling_vector_entities}
	\begin{array}{c}
		w_{2i-1} = w^{+}_{k,m},\quad 
		w_{2i} = w^{-}_{k,m},\quad
		i = (M+1)(k+N-1)+m+1.
	\end{array}
\end{equation}

\subsection{One-dimensional Variational Problem}\label{sub:402}

In the previous section, the traveling waves $w^{\pm}_{k,m}$ have been linearly expressed  through the free vector-valued function $\boldsymbol{y}$ as well as the initial and terminal states $\big(v_{0}(x),r_{0}(x)\big)$ and $\big(v_{1}(x),r_{1}(x)\big)$. 
It can be proven \cite{Kostin:2022} that we can symbolically represent this relation according to
\begin{equation}\label{eq:402:linear_state_relation}
	\boldsymbol{w}(\boldsymbol{y}(z),z,c_{1}) = \boldsymbol{A}\boldsymbol{y}(z) + c_{1}\boldsymbol{a} +\boldsymbol{g}(z),\quad
	z\in[0,\lambda]. 
\end{equation}		
Here,   $\boldsymbol{A} \in \mathbb{R}^{N_{w}\times N_{s}}$ is a known constant matrix, $\boldsymbol{a} \in \mathbb{R}^{N_{w}}$ is a known vector, and $\boldsymbol{g}:[0,\lambda]\to \mathbb{R}^{N_{w}}$ is a function defined by the initial and terminal values of $v$ and $r$. 	

Let us also take into account the vertex conditions in~\eqref{eq:310:odd_links}, \eqref{eq:310:even_links} and represent them as
\begin{equation}\label{eq:402:essential_boundary_conditions}
	\begin{array}{c}
		\boldsymbol{B}_{1}\boldsymbol{y}(\lambda) - \boldsymbol{B}_{0}\boldsymbol{y}(0) = 
		c_{1}\boldsymbol{b}_{1} + \boldsymbol{b}_{0},
		\\[1ex]
		\{\boldsymbol{B}_{0},\boldsymbol{B}_{1}\} \subset \mathbb{R}^{N_{b}\times N_{s}},\quad
		\{\boldsymbol{b}_{0},\boldsymbol{b}_{1}\} \subset \mathbb{R}^{N_{b}}. 
	\end{array}
\end{equation}			
In these terms, the optimal control problem~\eqref{eq:209:energy_minimization} can be reduced to a problem of calculus of variations:
Find such a function $\boldsymbol{y}^{*}(z)$  over $z\in[0,\lambda]$ and a constant $c^{*}_{1}$ that minimize the functional
\begin{equation} \displaystyle
E[\boldsymbol{y}^{*}] = \min\limits_{\boldsymbol{y}, c_{1}} E[\boldsymbol{y}]
\end{equation}
subject to the boundary constraints~\eqref{eq:402:essential_boundary_conditions}.

It can be derived from~\eqref{eq:401:term_mean_energy} and~\eqref{eq:402:energy_decomposition} that  
that the Euler--Lagrange ODE system is given by
\begin{equation}\label{eq:403:Euler_Lagrange_ODE}
	\begin{array}{c}
		\boldsymbol{A}^{\mathrm{T}}\boldsymbol{A}\boldsymbol{y}^{\prime\prime}(z) =
		-\boldsymbol{A}^{\mathrm{T}}\boldsymbol{g}^{\prime\prime}(z),\quad
		z\in[0,\lambda].
	\end{array}
\end{equation}			
The conjugate vector-valued function to the vector $\boldsymbol{y}$ is given as 
\begin{equation}
		\boldsymbol{p} = \partial e/\partial \boldsymbol{y}^{\prime} = \boldsymbol{A}^{\mathrm{T}}\boldsymbol{A}\boldsymbol{y}^{\prime} + \boldsymbol{A}^{\mathrm{T}}\boldsymbol{g}^{\prime}.
\end{equation}
The natural boundary condition in variations is
\begin{equation}\label{eq:404:boundary_transversality_conditions}
	\begin{array}{c}
		\boldsymbol{p}(\lambda )\cdot \delta\boldsymbol{y}(\lambda) - \boldsymbol{p}(0)\cdot \delta\boldsymbol{y}(0)   = 0.
	\end{array}
\end{equation}			
Next, we take into account the variation of essential boundary constraints~\eqref{eq:402:essential_boundary_conditions} according to
\begin{equation}\label{eq:404:variation_boundary_conditions}
	\boldsymbol{B}_{1}\delta\boldsymbol{y}(\lambda) - \boldsymbol{B}_{0}\delta\boldsymbol{y}(0) = \delta c_{1}\boldsymbol{b}_{1}.
\end{equation}			
From~\eqref{eq:404:variation_boundary_conditions}, we can express the variation $\delta c_{1}$ and represent
the natural conditions on $\boldsymbol{p}(0)$ and $\boldsymbol{p}(\lambda)$ 
in the form
\begin{equation}\label{eq:404:natural_boundary_conditions}
	\begin{array}{c}
		\boldsymbol{p}(0) = \boldsymbol{C}_{0}^{\mathrm{T}}\boldsymbol{h},\quad
		\boldsymbol{p}(\lambda)  = \boldsymbol{C}_{1}^{\mathrm{T}}\boldsymbol{h},
		\\ 
		\boldsymbol{C}_{i} = \boldsymbol{B}_{i} - \
		|\boldsymbol{b}_{1}|^{-2}\boldsymbol{b}\boldsymbol{b}^{\mathrm{T}}_{1}
		\boldsymbol{B}_{i},\quad
		i=0,1 ,\quad
		\boldsymbol{h} \in \mathbb{R}^{N_{b}}.
	\end{array}
\end{equation}				
Finally, the system of the boundary constraints~\eqref{eq:402:essential_boundary_conditions}, \eqref{eq:404:natural_boundary_conditions} is solved w.r.t. $2N_{s}$ constants appearing after integration of~\eqref{eq:403:Euler_Lagrange_ODE} and $N_{b}$ unknown components of the vector $\boldsymbol{h}$.

\subsection{Example of Optimal Control Design}\label{sub:405}

\begin{figure}[t]
	\begin{center}
		\includegraphics[width=0.6\linewidth]{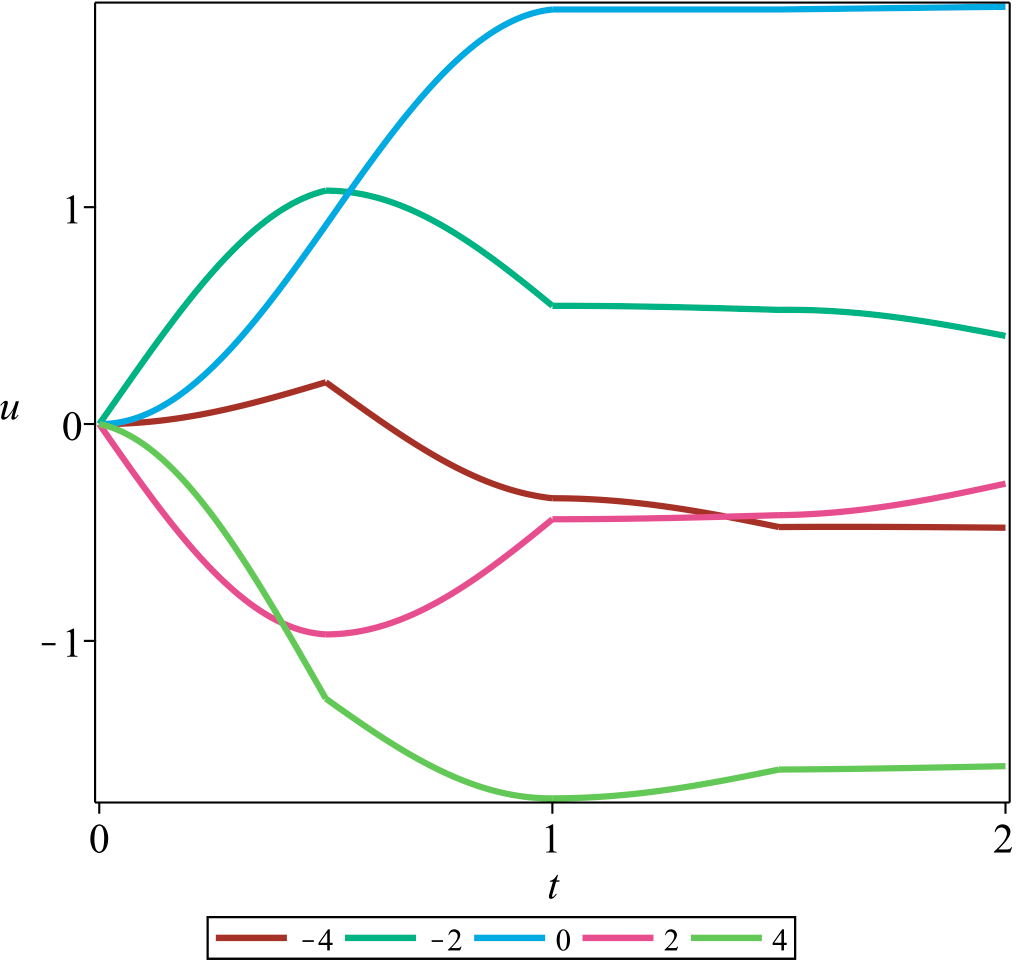}\\ 
		\caption{Optimal control $u_{n}(t)$ for $n\in J_{x}$ on the dimensionless time-space domain $\Omega = (0, 2)\times(-1,1)$ with  $N=4$.}
		\label{fig:10}
	\end{center}
\end{figure}

\begin{figure}[t]
	\begin{center}
		\includegraphics[width=0.6\linewidth]{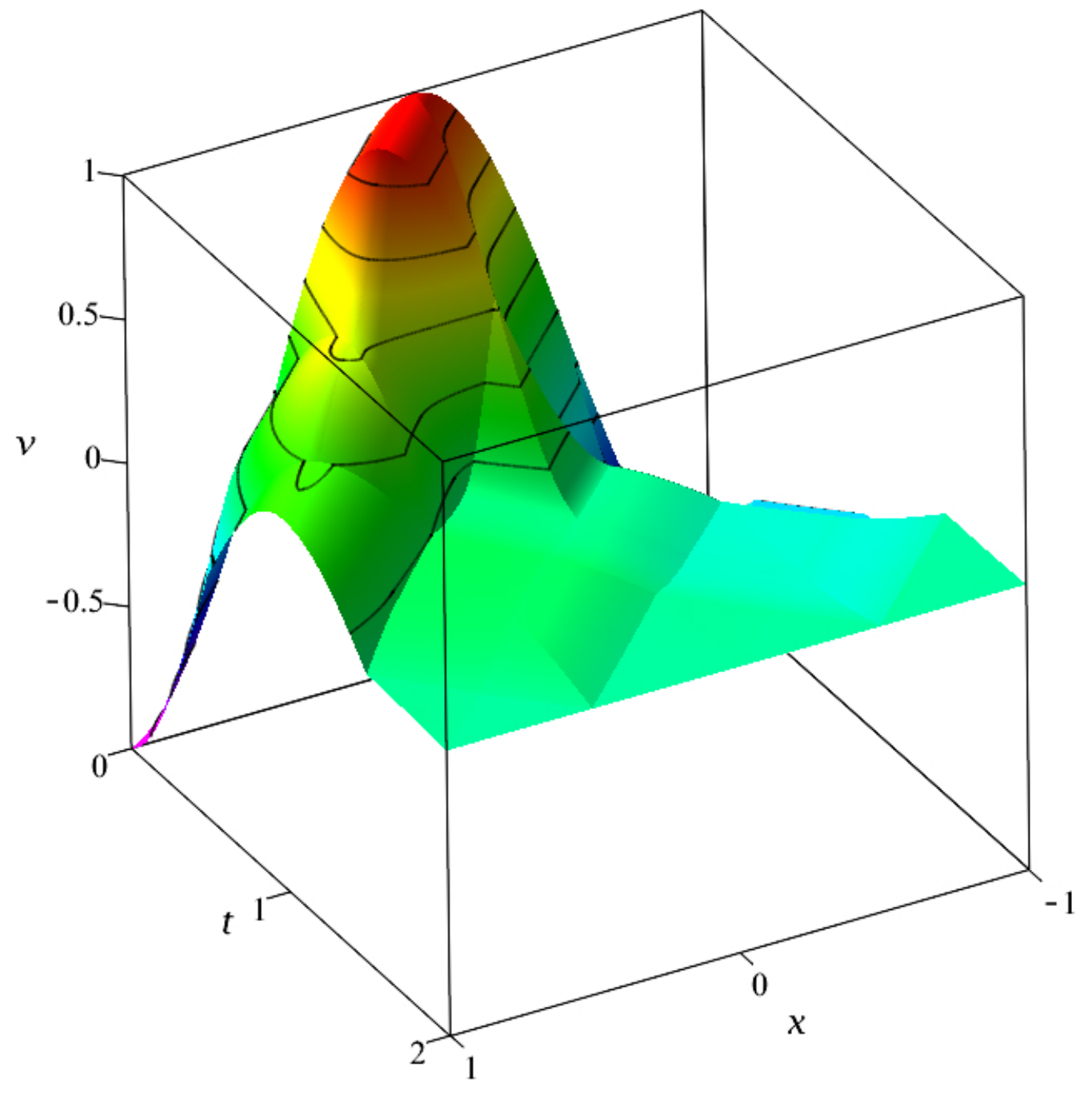}\\ 
		\caption{Optimal displacements $v(t,x)$ at $N=4$ and $T=2$.}
		\label{fig:08}
	\end{center}
\end{figure}

Let us consider an example of optimal control designed for the system with four distributed inputs ($N=4$).
The control time is equal to $T = M\lambda = 2$, what generates the mesh presented  in Fig.~\ref{fig:02}. The mesh parameters are given in Table~\ref{tb:parameters}.
\begin{table}[hb]
\begin{center}
\caption{Dimensionless parameters}\label{tb:parameters}
\begin{tabular}{ccccc}
$T$ (time) & $N$ & $M$ & $\lambda$ (length) \\\hline
$2$ &  $4$ & $4$ & $1/2$ \\ \hline
\end{tabular}
\end{center}
\end{table}
The initial and terminal conditions are taken as 
\begin{equation}\label{eq:405:initial_conditions} 
	\begin{array}{c}
		v(0,x)=\cos 3x,\quad
		r(0,x)=-\cos 3x,
	\quad
		v(T,x)=0,\quad
		r(T,x)=c_{1}.
	\end{array}
\end{equation}
The optimal displacements $v(t,x)$ are shown in Fig.~\ref{fig:08}. As seen from the graph, the rod arrives to the zero state at the terminal time instant $t=2$ both with zero energy and position. The corner points of the displacement function are visible, located according to the characteristics, which are shown in in Fig.~\ref{fig:02} by dashed lines.

The optimal integrals of the force jumps $u_{n}(t)$ with $n\in J_{x} = \{-4,-2,0,2,4\}$, which are components of the control vector $\boldsymbol{u}(t)$, are shown in Fig.~\ref{fig:10}.
According to the general solution of  Euler--Lagrange equations and the chosen initial and  terminal conditions, these integrals are linear combinations of trigonometric  and polynomial functions of time. Each control map $u_{n}(t)$ is continuous and by definition in~\eqref{eq:206:control_integrals} has got the zero initial value $u_{n}(0)=0$.

To restore the potential $r(t,x)$ and the external loads $f_{\pm}(t)$ as well as internal forces $f(t,x)$, we need first to calculate the control integrals $u_{k}(t)$ for $k\in J_{c}=\{-5,-3,-1,1,3,5\}$. 
These functions are obtained from the optimal integrals $u_{n}(t)$ with $n\in J_{x}$. 
The control inputs 
are found as the first derivatives of control integrals $f_{k}(t)=u'_{k}(t)$ with $k\in J_{c}$. 
These functions have discontinuity at the three time points   
$t=\frac12,1,\frac32$.

Let us analyze the change in the optimal value of the 
cost functional $E$ depending on the number of elements $N$ and control time $T$. Since the minimized mean energy $E$ introduced in~\eqref{eq:209:energy_minimization} is the integral of the energy density $e(t,x)$ over the domain $\Omega$ divided by the fixed parameter $T$, it is worthwhile 
to look into dependence of the minimum of this integral $T\cdot E$ on the control parameters. In the case of the discrete time set $T = M\lambda$, we can use the mesh numbers $M$ and $N$ as independent parameters.

\begin{figure}[t]
	\begin{center}
		\includegraphics[width=0.6\linewidth]{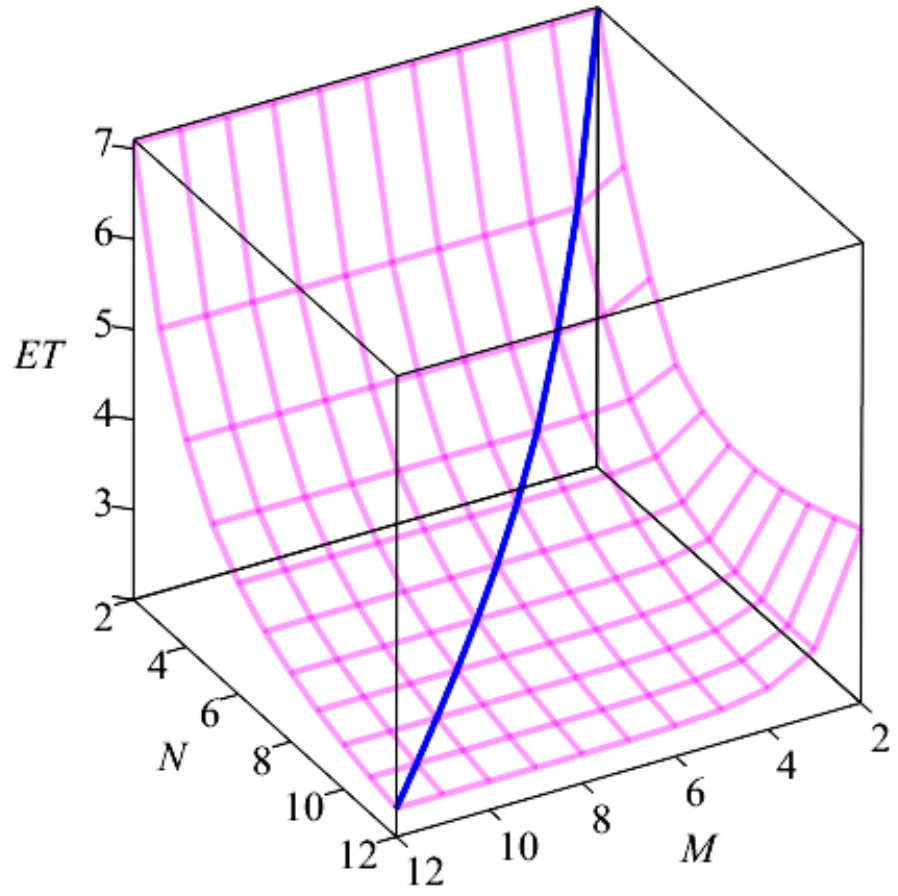}\\ 
		\caption{Energy integral $TE$ vs. the mesh numbers $M$, $N$.}
		\label{fig:14}
	\end{center}
\end{figure}

The optimal values of the energy integral $T\cdot E$ versus the mesh parameters $M$, $N$ are depicted in Fig.~\ref{fig:14} for the initial and terminal distribution~\eqref{eq:405:initial_conditions}. As seen on the graph, the functional $T\cdot E(M,N)$ is monotonically decreasing with increasing of both the number of time intervals $M$ and the number of space segments $N$. After a few time steps, the  rate of change of the functional w.r.t. $M$ becomes almost negligible. In contrast, the value $T\cdot E$ continues to decline with an increase in the number of elements $N$.  

To demonstrate the behavior of the function $T\cdot E(M,N)$ at some control time $T=\mathrm{const}$, we may fix the ratio $\frac{M}{N}$.  For example, the isochronic value $\frac{M}{N}=1$ gives us the control time 
$T=2$. 
The blue line in Fig.~\ref{fig:14} presents a monotonic decline of the one-dimensional function $T\cdot E(N,N)$. 
The case $N=2$ is special because the optimal value of the considered functional is the same for all $M \ge 2$.  
It is also worth noting the asymptotic behavior of the mapping $T\cdot E(2,N)$, which seems to tend to the limiting value. The proof of this fact will require a correct transition to the limiting case $N=\infty$ of distributed control, which is beyond the scope of this work.

\section{Conclusions}\label{sec:5}

The motion of a dynamic system under external boundary loads and internal distributed stresses has been studied. 
The mathematical model proposed can be related to longitudinal vibrations of a thin elastic rod controlled by forces applied along its central line together with normal forces acting at the ends. 
For given initial and terminal states and a fixed time horizon, the optimal control problem is to minimize the mean energy stored in the rod during the motion.  
In the case of equidistantly placed actuators and a uniform rod, the shortest possible time for bringing the system with any number of control elements to an arbitrary state is defined. An optimization algorithm using traveling waves is proposed to reduce the original problem to an one-dimensional variational problem with  boundary conditions of a special kind. The latter is solved exactly and the optimal control inputs explicitly.
A practical implementation of such a feedforward control usually requires additional feedback control to ensure robustness. We plan to enhance our approach by taking into account that piezoelectric elements, which can be utilized for  realization, serve both as sensors and observers, and use that information for additional feedback stabilization.

\end{document}